\documentclass[12pt]{amsart}
\usepackage{amssymb}

\usepackage[]{amsfonts}

\def\underset#1#2{\mathrel{\mathop{\kern0pt #2}\limits_{#1}}}

\def\overset#1#2{\mathrel{\mathop{\kern0pt #2}\limits^{#1}}}

\def\couleur(#1 #2 #3)
	{
\special{ps::[end]}}

\def\sqr#1#2{{\vcenter{\vbox{\hrule height.#2pt
             \hbox{\vrule width.#2pt height#1pt \kern#1pt
             \vrule width.#2pt}
             \hrule height.#2pt}}}}

\def\square{\mathchoice\sqr64\sqr64\sqr43\sqr23}			

\def\st{\mathinner{\mkern1mu\raise1pt\hbox{.}				
		   \mkern1mu\raise4pt\hbox{.}
		   \mkern1mu\raise1pt\hbox{.}
		 }
         }

\def\bx#1{\setbox1=\hbox{\kern3pt{#1}\kern3pt}				
 \dimen1=\ht1 \advance\dimen1 by 3pt \dimen2=\dp1 \advance\dimen2 by 3pt
 \setbox1=\hbox{\vrule height\dimen1 depth\dimen2\box1\vrule}%
 \setbox1=\vbox{\hrule\box1\hrule}%
 \advance\dimen1 by .4pt \ht1=\dimen1
 \advance\dimen2 by .4pt \dp1=\dimen2 \box1\relax}

\def\k#1{\kern#1em}
\def\vci{\vrule  width.02em height1.47ex depth-.0ex}				
\def\11{{\rm\k{.2}\vci\k{-.37}1}}

\newtheorem{Theorem}{Theorem}[section]
\newtheorem{Definition}[Theorem]{Definition}
\newtheorem{Remark}[Theorem]{Remark}
\newtheorem{Corollary}[Theorem]{Corollary}
\newtheorem{Lemma}[Theorem]{Lemma}

\begin{document}
\title{On linear extension for interpolating sequences.}
\author{Eric Amar}
\address{University Bordeaux I, 351 Cours de la Lib{\'e}ration, 33405 
Talence France}
\email{Eric.Amar@math.u-bordeaux1.fr{\hskip 1.8em}}
\maketitle
\begin{abstract} {
Let $A$ be a uniform algebra on the compact space $X$ and $\sigma $ a 
probability measure on $X$. We define the Hardy spaces $H^{p}(\sigma )$ 
and the $\displaystyle H^{p}(\sigma )$ interpolating sequences $S$ in 
the $p$-spectrum ${\mathcal{M}}_{p}$ of $\sigma $. We prove, under some 
structural hypotheses on $\sigma $ that "Carleson type" conditions on 
$S$ imply that $S$ is interpolating with a linear extension operator in 
$\displaystyle H^{s}(\sigma ),\ s<p$ provided that either $p=\infty $ 
or $p\leq 2$.\ \par
{\hskip 1.8em}This gives new results on interpolating sequences for Hardy 
spaces of the ball and the polydisc. In particular in the case of the 
unit ball of ${\mathbb{C}}^{n}$ we get that if there is a sequence $\{\rho _{a}\}_{a\in S}$ 
bounded in $\displaystyle H^{\infty }({\mathbb{B}})$ such that $\forall a,b\in S,\ \rho _{a}(b)=\delta _{ab}$, 
then $S$ is $\displaystyle H^{p}({\mathbb{B}})$-interpolating with a linear 
extension operator for any $1\leq p<\infty $.\ \par
}\end{abstract}
\section{Introduction{\hskip 1.8em}}
\setcounter{equation}{0}Let $\displaystyle {\mathbb{B}}$ be the unit ball 
of ${\mathbb{C}}^{n}$; we denote as usual by $\displaystyle H^{p}({\mathbb{B}})$ 
the Hardy spaces of holomorphic functions in $\displaystyle {\mathbb{B}}$. 
Let $S$ a sequence of points in $\displaystyle {\mathbb{B}}$ and $1\leq p\leq \infty $ 
; we say that $S$ is $H^{p}$-interpolating if\ \par

\begin{displaymath} 
\forall \lambda \in \ell ^{p}(S),\ \exists f\in H^{p}({\mathbb{B}})\ s.t.\ 
\forall a\in S,\ f(a)=\lambda _{a}(1-\displaystyle \left\vert{a}\right\vert 
^{2})^{n/p}.\end{displaymath} \ \par
{\hskip 1.8em}Let $a\in {\mathbb{B}}$ we set $\displaystyle k_{a}(z):=\displaystyle \frac{1}{(1-\overline{a}\cdot z)^{n}}$ 
its reproducing kernel and $\displaystyle k_{p,a}:=\displaystyle \frac{k_{a}}{\displaystyle \left\Vert{k_{a}}\right\Vert 
_{p}}$ the normalized reproducing kernel for $a$ in $\displaystyle H^{p}({\mathbb{B}})$. 
Now if $S$ is $H^{p}$-interpolating, then we have, with $p'$ the conjugate 
exponent for $p$:\ \par
{\hskip 3.6em}$\displaystyle \exists C>0,\ \forall a\in S,\ \exists \rho _{a}\in H^{p}({\mathbb{B}})\ 
s.t.\ \displaystyle \left\langle{\rho _{a},\ k_{p',b}}\right\rangle =\delta 
_{ab}$.\ \par
{\hskip 1.8em}We shall say that $S$ is dual bounded in $\displaystyle 
H^{p}({\mathbb{B}})$ if the dual system $\{\rho _{a}\}_{a\in S}$ to $\{k_{p',a}\}_{a\in S}$ 
exits and is bounded in $\displaystyle H^{p}({\mathbb{B}})$.\ \par
Hence if $S$ is $H^{p}$-interpolating then $S$ is dual bounded in $\displaystyle 
H^{p}({\mathbb{B}})$.\ \par
\begin{Definition} We say that the $\displaystyle H^{p}({\mathbb{B}})$ 
interpolating sequence $S$ has the linear extension property (LEP) if 
there is a bounded linear operator $E\ :\ \ell ^{p}{\longrightarrow}H^{p}({\mathbb{B}})$ 
such that $\displaystyle \forall \lambda \in \ell ^{p},\ E\lambda $ interpolates 
the sequence $\lambda $ in $\displaystyle H^{p}({\mathbb{B}})$ on $S$, 
i.e.\ \par
{\hskip 3.6em}$\displaystyle \exists C>0,\ \forall \lambda \in \ell ^{p},\ E\lambda \in H^{p}({\mathbb{B}}),\ 
\displaystyle \left\Vert{E\lambda }\right\Vert _{p}\leq C\ s.t.\ \forall 
a\in S,\ E\lambda (a)=\lambda _{a}\displaystyle \left\Vert{k_{a}}\right\Vert 
_{p'}$\ \par
\end{Definition}
\ \par
{\hskip 1.8em}Natural questions are the following:\ \par
{\hskip 1.8em}If $S$ is dual bounded in $\displaystyle H^{p}({\mathbb{B}})$, 
is $S\in IH^{p}({\mathbb{B}})$ ?\ \par
{\hskip 1.8em}If $S\in IH^{p}({\mathbb{B}})$ has $S$ automatically the 
LEP ?\ \par
\ \par
{\hskip 1.8em}This is true in the classical case of the Hardy spaces of 
the unit disc ${\mathbb{D}}$ :\ \par
{\hskip 1.8em}for $p=\infty $ this is the famous characterization of $\displaystyle 
H^{\infty }$ interpolating sequences by L. Carleson~\cite{CarlInt58} and 
the LEP was given by P. Beurling~\cite{PBeurling62}.\ \par
{\hskip 1.8em}for $p\in [1,\ \infty [$ this was done by H. Shapiro and 
A. Shields~\cite{ShapShields61} and because the characterization is the 
same for all $p\in [1,\infty ]$, the LEP is deduced easily from the $\displaystyle 
H^{\infty }$ case and was done explicitly with $\overline{\partial }$ 
methods in~\cite{amExt83}.\ \par
{\hskip 1.8em}For the Bergman classes $A^{p}({\mathbb{D}})$, it is no 
longer true that the interpolating sequences are the same for $A^{p}({\mathbb{D}})$ 
and $A^{q}({\mathbb{D}}),\ q{\not =}p$. But A.P. Schuster and K. Seip~\cite{seip2},~\cite{SchusterSeip00} 
proved that $S$ dual bounded in $A^{p}({\mathbb{D}})$ implies $S\ A^{p}({\mathbb{D}})$-interpolating 
still with the LEP.\ \par
{\hskip 1.8em}The first question is open, even in the ball ${\mathbb{B}}$ 
of ${\mathbb{C}}^{n},\ n\geq 2$, with $H^{p}({\mathbb{B}})$, the usual 
Hardy spaces of the ball or in the polydisc ${\mathbb{D}}^{n}$ of ${\mathbb{C}}^{n},\ n\geq 2$ 
still with the usual Hardy spaces.\ \par
{\hskip 1.8em}The second one is known only in the case $p=\infty $ as 
we shall see later.\ \par
{\hskip 1.8em}Nevertheless in the case of the unit ball of ${\mathbb{C}}^{n}$, 
B. Berndtsson~\cite{Bernd85} proved that if the product of the Gleason 
distances of the points of $S$ is bounded below away of $0$ then $S$ is 
$\displaystyle H^{\infty }({\mathbb{B}})$. He also proved that this condition 
is not necessary for $n>1$.\ \par
{\hskip 1.8em}B. Berndtsson, A. S-Y. Chang and K-C. Lin~\cite{BernChanLin87} 
proved the same theorem in the polydisc of ${\mathbb{C}}^{n}$.\ \par
{\hskip 1.8em}In this paper we shall prove that loosing a little bit on 
the value of $p$, $S$ dual bounded in $\displaystyle H^{p}({\mathbb{B}})$ 
implies $\forall s<p,\ S\in IH^{s}({\mathbb{B}})$ with the LEP, provided 
that $1<p\leq 2$ or $p=\infty $. In particular:\ \par
\begin{Theorem} If $S\subset {\mathbb{B}}$ is dual bounded in $\displaystyle 
H^{p}({\mathbb{B}})$, then it is $\displaystyle H^{s}$-interpolating for 
any $1\leq s<p$, provided that $p\in ]1,2]$ or $p=\infty $. Moreover $S$ 
has the property that there is a bounded linear operator from $\ell ^{s}(S){\longrightarrow}H^{s}({\mathbb{B}})$ 
doing the interpolation.\ \par
\end{Theorem}
{\hskip 1.8em}The methods we use being purely functional analytic, these 
results extend to the setting of uniform algebras.\ \par
\section{Uniform algebras.{\hskip 1.8em}}
\setcounter{equation}{0}Let $A$ be a uniform algebra on the compact space 
$X$, i.e. $A$ is a sub-algebra of ${\mathcal{C}}(X)$, the continuous functions 
on $X$, which separates the points of $X$ and contains $1$.\ \par
{\hskip 1.8em}Let $\sigma $ be a probability measure on $X$.\ \par
{\hskip 1.8em}For $1\leq p<\infty $ we define as usual the Hardy space 
$\displaystyle H^{p}(\sigma )$ as the closure of $A$ in $L^{p}(\sigma )$.\ 
\par
$\displaystyle H^{\infty }(\sigma )$ will be the weak-* closure of $A$ 
in $L^{\infty }(\sigma )$.\ \par
{\hskip 1.8em}Let ${\mathcal{M}}$ be the Guelfand spectrum of $A$, i.e. 
the multiplicative elements of $A'$. We note the same way an element of 
$A$ and its Guelfand transform:\ \par
{\hskip 3.9em}$\forall a\in {\mathcal{M}}\subset A',\ \forall f\in A,\ f(a):=\hat f(a)=a(f)$.\ 
\par
{\hskip 1.8em}We shall use the following notions, already introduced in~\cite{amIntInt06}.\ 
\par
\begin{Definition} Let ${\mathcal{M}}$ be the spectrum of $A$ and $a\in {\mathcal{M}}$; 
we call $k_{a}\in H^{p}(\sigma )$ a $p$-reproducing kernel for the point 
$a$ if $\displaystyle \forall f\in A,\ f(a)=\displaystyle \int_{X}^{}{f(\zeta )\overline{k}_{a}(\zeta 
)\,d\sigma (\zeta )}$.\ \par
{\hskip 1.8em}We define the $p$-spectrum of $\sigma $ as the subset $\displaystyle 
{\mathcal{M}}_{p}$ of ${\mathcal{M}}$ such that every element has a $p'$-reproducing 
kernel with $p'$ the conjuguate exponent for $p$, $\displaystyle \displaystyle \frac{1}{p}+\displaystyle \frac{1}{p'}=1$.\ 
\par
\end{Definition}
\begin{Definition} We say that $S\subset {\mathcal{M}}_{p}$ is $\displaystyle 
H^{p}(\sigma )$ interpolating for $1\leq p<\infty $, $S\in IH^{p}(\sigma )$ 
if\ \par
{\hskip 3.6em}$\displaystyle \forall \lambda \in \ell ^{p},\ \exists f\in H^{p}(\sigma )\ s.t.\ \forall 
a\in S,\ f(a)=\lambda _{a}\displaystyle \left\Vert{k_{a}}\right\Vert _{p'}$.\ 
\par
We say that $S\subset {\mathcal{M}}_{\infty }$ is $\displaystyle H^{\infty }(\sigma )$ 
interpolating, $S\in IH^{\infty }(\sigma )$ if\ \par
{\hskip 3.6em}$\displaystyle \forall \lambda \in \ell ^{\infty },\ \exists f\in H^{\infty }(\sigma 
)\ s.t.\ \forall a\in S,\ f(a)=\lambda _{a}$.\ \par
\end{Definition}
\begin{Remark} If $S$ is $\displaystyle H^{p}(\sigma )$-interpolating 
then there is a constant $C_{I}$, the interpolating constant, such that~\cite{amIntInt06}:\ 
\par
{\hskip 1.8em}$\displaystyle \forall \lambda \in \ell ^{p},\ \exists f\in H^{p}(\sigma ),\ \displaystyle 
\left\Vert{f}\right\Vert _{p}\leq C_{I}\displaystyle \left\Vert{\lambda 
}\right\Vert _{p},\ s.t.\ \forall a\in S,\ f(a)=\lambda _{a}\displaystyle 
\left\Vert{k_{a}}\right\Vert _{p'}.$\ \par
\end{Remark}
\begin{Definition} We say that the $\displaystyle H^{p}(\sigma )$ interpolating 
sequence $S$ has the linear extension property (LEP) if there is a bounded 
linear operator $E\ :\ \ell ^{p}{\longrightarrow}H^{p}(\sigma )$ such 
that $\displaystyle \forall \lambda \in \ell ^{p},\ E\lambda $ interpolates 
the sequence $\lambda $ in $\displaystyle H^{p}(\sigma )$ on $S$, i.e.\ 
\par
{\hskip 3.6em}$\displaystyle \exists C>0,\ \forall \lambda \in \ell ^{p},\ E\lambda \in H^{p}(\sigma 
),\ \displaystyle \left\Vert{E\lambda }\right\Vert _{p}\leq C\ s.t.\ \forall 
a\in S,\ E\lambda (a)=\lambda _{a}\displaystyle \left\Vert{k_{a}}\right\Vert 
_{p'}$\ \par
\end{Definition}
{\hskip 1.8em}Let $S\subset {\mathcal{M}}_{p}$, so $\displaystyle k_{p',a}:=\displaystyle \frac{k_{a}}{\displaystyle \left\Vert{k_{a}}\right\Vert 
_{p'}}$, the normalized reproducing kernel, exits for any $a\in S$; let 
us consider a dual system $\{\rho _{a}\}_{a\in S}\subset H^{p}(\sigma )$, 
i.e. $\displaystyle \forall a,b\in S,\ \displaystyle \left\langle{\rho _{a},\ k_{p',b}}\right\rangle 
=\delta _{a,b}$ when it exists.\ \par
\begin{Definition} We say that $S\subset {\mathcal{M}}_{p}$ is dual bounded 
in $\displaystyle H^{p}(\sigma )$ if a dual system $\{\rho _{a}\}_{a\in S}\subset H^{p}(\sigma )$ 
exists and if this sequence is bounded in $\displaystyle H^{p}(\sigma )$, 
i.e. $\exists C>0\ s.t.\ \forall a\in S,\ \left\Vert{\rho _{a}}\right\Vert _{p}\leq 
C$.\ \par
\end{Definition}
{\hskip 1.8em}We shall show that, under some structural hypotheses on 
$\sigma $ and the fact that $S$ is Carleson (the definition of Carleson 
sequences will be given later):\ \par
\begin{Theorem}  \label{algUniExt613}If $1\leq s<p$ and either $p\leq 2$ 
or $p=\infty $, $S\subset {\mathcal{M}}_{p}\cap {\mathcal{M}}_{s}$ is 
dual bounded  in $\displaystyle H^{p}(\sigma )$ and $S$ is a Carleson 
sequence, then $S\in IH^{s}(\sigma )$ with the linear extension property.\ 
\par
\end{Theorem}
{\hskip 1.8em}The passage from $p=2$ to $p\leq 2$ in the case of the ball 
is due to F. Bayart: he uses Khintchine's inequalities which reveal to 
be very well fitted to this problem. In fact F. Lust-Piquart showed me 
a way not to use Khintchine's inequalities: one can use the fact that 
$L^{p}$ spaces are of type $p$ in the part $p\leq 2$ in the proof of theorem~\ref{algUniExt613}.\ 
\par
{\hskip 1.8em}I shall add this proof.\ \par
{\hskip 1.8em}The case $p=\infty $ of this theorem is the best possible 
in this generality. There is no hope to have that dual boundedness in 
$\displaystyle H^{\infty }$ implies $\displaystyle H^{\infty }$-interpolation 
as L. Carleson proved for the unit disc.\ \par
{\hskip 1.8em}In~\cite{GorMorLin90} and in~\cite{Izuchi91} the authors 
proved that in the spectrum of the uniform algebra $H^{\infty }({\mathbb{D}})$ 
there are sequences $S$ of points such that the product of the Gleason 
distances is bounded below away from $0$, which implies that $S$ is dual 
bounded in $H^{\infty }({\mathbb{D}})$, but $S$ is \textsl{not} $H^{\infty }$-interpolating.\ 
\par
\ \par
{\hskip 1.8em}The general theorem~\ref{algUniExt613} implies a polydisc 
and a ball version.\ \par
In the polydisc ${\mathbb{D}}^{n}\subset {\mathbb{C}}^{n}$ the structural 
hypotheses are true~\cite{amIntInt06}, hence\ \par
\begin{Theorem} Let $S\subset {\mathbb{D}}^{n}$ be a Carleson sequence 
and dual bounded  in $\displaystyle H^{p}({\mathbb{D}}^{n})$ with either 
$p=\infty $ or $1<p\leq 2$, then $S$ is $\displaystyle H^{s}({\mathbb{D}}^{n})$ 
interpolating for any $1\leq s<p$ with the LEP.\ \par
\end{Theorem}
{\hskip 1.8em}In the ball, the structural hypothesesare true~\cite{amIntInt06} 
and moreover we know, by an easy corollary of a theorem of P. Thomas~\cite{Thomas87}, 
that $S$ dual bounded in $\displaystyle H^{p}({\mathbb{B}})$ implies $S$ 
Carleson, hence\ \par
\begin{Theorem} \label{algUniInt3}Let $S\subset {\mathbb{B}}$ be dual 
bounded  in $\displaystyle H^{p}({\mathbb{B}})$ with either $p=\infty $ 
or $1<p\leq 2$, then $S$ is $\displaystyle H^{s}({\mathbb{B}})$ interpolating 
for any $1\leq s<p$ with the LEP.\ \par
\end{Theorem}
{\hskip 1.8em}As usual by use of the "subordination lemma"~\cite{amBerg78} 
we have the same result for the Bergman classes of the ball. Denote by 
$A^{p}({\mathbb{B}})$ the holomorphic functions in $L^{p}({\mathbb{B}})$ 
for the area measure of the ball then\ \par
\begin{Corollary} Let $S\subset {\mathbb{B}}$ be dual bounded  in $\displaystyle 
A^{p}({\mathbb{B}})$ with either $p=\infty $ or $1<p\leq 2$, then $S$ 
is $\displaystyle A^{s}({\mathbb{B}})$ interpolating for any $1\leq s<p$ 
with the LEP.\ \par
\end{Corollary}
{\hskip 1.8em}In~\cite{amIntInt06} it was proved:\ \par
\begin{Theorem} \label{algUniExt820}Let $p\geq 1$, $1\leq s<p$ and $q$ 
be such that $\displaystyle \displaystyle \frac{1}{s}=\displaystyle \frac{1}{p}+\displaystyle \frac{1}{q}$. 
Suppose that $S\subset {\mathcal{M}}_{s}\cap {\mathcal{M}}_{q'}$ is $\displaystyle 
H^{p}(\sigma )$ interpolating, $q$-Carleson and $\sigma $ verifies the 
structural hypotheses , then $S$ is $\displaystyle H^{s}(\sigma )$ interpolating.\ 
\par
\end{Theorem}
{\hskip 1.8em}The theorem~\ref{algUniExt820} is better for $p\in [1,\ 2]$ 
or $p=\infty $ : we have the LEP under the weaker assumption that $S$ 
is dual bounded in $\displaystyle H^{p}(\sigma )$.\ \par
{\hskip 1.8em}But we have not the full range of $p$ as in theorem~\ref{algUniExt820}.\ 
\par
\section{Reproducing kernels.{\hskip 1.8em}}
\setcounter{equation}{0}Let us recall some facts about reproducing kernels 
and $p$-spectrum.\ \par
{\hskip 1.8em}First the reproducing kernel for $a\in {\mathcal{M}}$ if 
it exists is unique. Suppose there are $2$ of them, $k_{a}\in H^{p}(\sigma )$ 
and $k'_{a}\in H^{q}(\sigma )$:\ \par
{\hskip 3.6em}$\displaystyle \forall f\in A,\ 0=f(a)-f(a)=\displaystyle \int_{X}^{}{f(\overline{k}_{a}-\overline{k}'_{a})\,ds}{\Longrightarrow}k_{a}=k'_{a}\ 
\sigma -a.e.$\ \par
because, by definition, $A$ is dense in $\displaystyle H^{r}(\sigma )$ 
with $r:=\min (p,q)$. Hence it is correct to denote it by $k_{a}$ without 
reference to the $\displaystyle H^{p}(\sigma )$ where it belongs.\ \par
{\hskip 1.8em}Let $a\in {\mathcal{M}}_{p}$ then $k_{a}\in H^{p'}(\sigma )$; 
if $p<q{\Longrightarrow}q'<p'$  hence $k_{a}\in H^{q'}(\sigma )$ because 
$\sigma $ is a probability measure so $a\in {\mathcal{M}}_{q}$ and we 
have $\displaystyle p<q{\Longrightarrow}{\mathcal{M}}_{p}\subset {\mathcal{M}}_{q}$.\ 
\par
{\hskip 1.8em}To simplify the notation we shall use:\ \par
{\hskip 3.6em}$\displaystyle \displaystyle \left\langle{f,g}\right\rangle :=\displaystyle \int_{X}^{}{f\overline{g}\,d\sigma 
},$\ \par
whenever this is meaningfull.\ \par
\ \par
{\hskip 1.8em}If $a\in {\mathcal{M}}_{2}$ we always have a "Poisson kernel" 
associated to $a$, $\displaystyle P_{a}:=\displaystyle \frac{\displaystyle \left\vert{k_{a}}\right\vert 
^{2}}{\displaystyle \left\Vert{k_{a}}\right\Vert _{2}^{2}}$ and the well 
known\ \par
\begin{Lemma} $P_{a}\in L^{1}(\sigma ),\ \left\Vert{P_{a}}\right\Vert _{1}=1$ 
and\ \par
{\hskip 3.6em}{\rm{}}$\displaystyle \forall f\in A,\ f(a)=\displaystyle \left\langle{f,\ P_{a}}\right\rangle 
=\displaystyle \int_{X}^{}{fP_{a}\,d\sigma }.$\ \par
\end{Lemma}
{\hskip 1.8em}Proof\ \par
{\hskip 1.8em}$\displaystyle \displaystyle \int_{X}^{}{fP_{a}\,d\sigma }=\displaystyle \int_{X}^{}{f\displaystyle 
\frac{k_{a}\overline{k}_{a}}{\displaystyle \left\Vert{k_{a}}\right\Vert 
_{2}^{2}}\,d\sigma }=\displaystyle \frac{1}{\displaystyle \left\Vert{k_{a}}\right\Vert 
_{2}^{2}}f(a)k_{a}(a)=f(a),$\ \par
because $fk_{a}\in H^{2}(\sigma )$ and $k_{a}(a)=\int_{X}^{}{k_{a}\overline{k}_{a}\,d\sigma }=\left\Vert{k_{a}}\right\Vert 
_{2}^{2}$. \hfill$\square$\ \par
{\hskip 1.8em}This allows us to define the Poisson integral of a bounded 
function on $X$:\ \par
\begin{Definition} Let $f\in L^{\infty }(\sigma )$ we set $\forall a\in {\mathcal{M}}_{2},\ \tilde f(a):=\left\langle{f,\ P_{a}}\right\rangle 
$ its Poisson integral.\ \par
If $f\in L^{2}(\sigma )$ we set $f^{*}:=P_{2}f$ its orthogonal projection 
on $\displaystyle H^{2}(\sigma )$; we extend $f^{*}$ on ${\mathcal{M}}_{2}$:\ 
\par
{\hskip 3.6em}$\displaystyle \forall f\in L^{2}(\sigma ),\ \forall a\in {\mathcal{M}}_{2},\ f^{*}(a):=\displaystyle 
\left\langle{f^{*},k_{a}}\right\rangle =\displaystyle \left\langle{f,\ 
k_{a}}\right\rangle .$\ \par
\end{Definition}
{\hskip 1.8em}Of course if $f\in A$ we have $f^{*}=\tilde f=f$ and for 
any $f\in L^{\infty }(\sigma ),\ \widetilde{(f^{*})}=f^{*}$.\ \par
\subsection{Structural hypotheses{\hskip 1.5em}}
We shall need some structural hypotheses on $\sigma $ relative to the 
reproducing kernels.\ \par
\begin{Definition} Let $q\in ]1,\infty [$, we say that the measure $\sigma $ 
verifies the structural hypothesis $SH(q)$ if, with $q'$ the conjugate 
of $q$:\ \par
{\hskip 1.8em}
\begin{equation} 
\exists \alpha =\alpha _{q}>0\ s.t.\ \forall a\in {\mathcal{M}}_{q}\cap 
{\mathcal{M}}_{q'}\subset {\mathcal{M}}_{2},\ \displaystyle \left\Vert{k_{a}}\right\Vert 
_{2}^{2}\geq \alpha \displaystyle \left\Vert{k_{a}}\right\Vert _{q}\displaystyle 
\left\Vert{k_{a}}\right\Vert _{q'}.
\end{equation} \ \par
\end{Definition}
This is opposite to the H{\"o}lder inequalities.\ \par
Because $a\in {\mathcal{M}}_{q}\cap {\mathcal{M}}_{q'}\subset {\mathcal{M}}_{2}$, 
we have $\displaystyle k_{a}(a)=\displaystyle \int_{X}^{}{k_{a}(\zeta )\overline{k}_{a}(\zeta 
)\,d\sigma }=\displaystyle \left\Vert{k_{a}}\right\Vert _{2}^{2}$ and 
the condition above is the same as\ \par
{\hskip 3.6em}$\left\Vert{k_{a}}\right\Vert _{q}\left\Vert{k_{a}}\right\Vert _{q'}\leq 
\alpha _{q}^{-1}k_{a}(a).$\ \par
\begin{Definition} Let $p,s\in [1,\infty ]$ and $q$ such that $\displaystyle 
\displaystyle \frac{1}{s}=\displaystyle \frac{1}{p}+\displaystyle \frac{1}{q}$. 
We say that the measure $\sigma $ verifies the structural hypothesis $SH(p,s)$ 
if\ \par
{\hskip 3.6em}{\rm{}}
\begin{equation} 
\exists \beta =\beta _{p,q}>0\ s.t.\ \forall a\in {\mathcal{M}}_{s},\ 
\displaystyle \left\Vert{k_{a}}\right\Vert _{s'}\leq \beta \displaystyle 
\left\Vert{k_{a}}\right\Vert _{p'}\displaystyle \left\Vert{k_{a}}\right\Vert 
_{q'}.
\end{equation} \ \par
\end{Definition}
This is meaningfull because $s<p,\ s<q$ hence ${\mathcal{M}}_{s}\subset {\mathcal{M}}_{p}\cap {\mathcal{M}}_{q}$.\ 
\par
{\hskip 1.8em}In the case of the unit ball ${\mathbb{B}}\subset {\mathbb{C}}^{n}$ 
and $\sigma $ the Lebesgue measure on $X=\partial {\mathbb{B}}$ and in 
the case of the polydisc ${\mathbb{D}}^{n}\subset {\mathbb{C}}^{n}$ and 
$\sigma $ the Lebesgue measure on ${\mathbb{T}}^{n}$, it is shown in~\cite{amIntInt06} 
that these two hypotheses are verified for all $p,s,q$.\ \par
\subsection{Interpolating sequences.{\hskip 1.8em}}
We shall use the following facts proved in~\cite{amIntInt06} :\ \par
\begin{Theorem} If, for a $p\geq 1$, $S\subset {\mathcal{M}}_{p}$, if 
$S\in IH^{\infty }(\sigma )$ and if $\sigma $ verifies $SH(p)$  then $S\in IH^{p}(\sigma )$ 
with the $L.E.P.$.\ \par
\end{Theorem}
\begin{Theorem} \label{algUniInt5}If $S\subset {\mathcal{M}}_{1}$ and 
$S$ is dual bounded in $\displaystyle H^{p}(\sigma )$ for a $p>1$, then 
$S\in IH^{1}(\sigma )$.\ \par
\end{Theorem}
{\hskip 1.8em}We shall need to truncate $S$ to its first $N$ elements, 
say $S_{N}$. Clearly if $S\in IH^{p}(\sigma )$ then $S_{N}\in IH^{p}(\sigma )$ 
with a smaller constant than $C_{I}$. Let $\displaystyle I_{S_{N}}^{p}:=\{f\in H^{p}(\sigma )\ s.t.\ f_{\mid S_{N}}=0\}$ 
be the module over $A$ of the functions zero on $S_{N}$. We have then 
for $\lambda \in \ell ^{p}$, with $\{\rho _{a}\}_{a\in S}$ a bounded dual 
sequence, that the function $f_{N}:=\sum_{a\in S_{N}}^{}{\lambda _{a}\rho _{a}}$ 
interpolates $\lambda $ on $S_{N}$ and we have $\left\Vert{f_{N}}\right\Vert _{H^{p}(\sigma )/I_{S_{N}}^{p}}\leq C_{I}\left\Vert{\lambda 
}\right\Vert _{p}$.\ \par
{\hskip 1.8em}We also have the converse for $1<p\leq \infty $, which is 
all what we need~\cite{amIntInt06}:\ \par
\begin{Lemma} \label{algUniInt0}If $S$ is such that all its truncations 
$S_{N}$ are in $IH^{p}(\sigma )$ for a $p>1$, with a uniform constant 
$C_{I}$ then $S\in IH^{p}(\sigma )$ with the same constant.\ \par
\end{Lemma}
\section{Carleson sequences.{\hskip 1.8em}}
\setcounter{equation}{0}As before we denote by $\displaystyle k_{q,a}:=\displaystyle \frac{k_{a}}{\displaystyle \left\Vert{k_{a}}\right\Vert 
_{q}}$ the normalized reproducing kernel in $\displaystyle H^{q}(\sigma )$.\ 
\par
\begin{Definition} We say that the sequence $S\subset {\mathcal{M}}_{q'}$ 
is a $q$-Carleson sequence if $1\leq q<\infty $ and\ \par
{\hskip 3.6em}$\displaystyle \exists D_{q}>0,\ \forall \mu \in \ell ^{q},\ \displaystyle \left\Vert{\displaystyle 
\sum_{a\in S}^{}{\mu _{a}k_{q,a}}}\right\Vert _{q}\leq D_{q}\displaystyle 
\left\Vert{\mu }\right\Vert _{q}.$\ \par
We say that the sequence $S\subset {\mathcal{M}}_{q'}$ is a weakly  $q$-Carleson 
sequence if $2\leq q<\infty $ and\ \par
{\hskip 3.6em}$\displaystyle \exists D_{q}>0,\ \forall \mu \in \ell ^{q},\ \displaystyle \left\Vert{\displaystyle 
\sum_{a\in S}^{}{\displaystyle \left\vert{\mu _{a}}\right\vert ^{2}\displaystyle 
\left\vert{k_{q,a}}\right\vert ^{2}}}\right\Vert _{q/2}\leq D_{q}\displaystyle 
\left\Vert{\mu }\right\Vert _{q}^{2}.$\ \par
\end{Definition}
{\hskip 1.8em}We call "weakly" Carleson the second condition because\ 
\par
\begin{Lemma} If $2\leq q<\infty $ and $S$ is $q$-Carleson then it is 
weakly $q$-Carleson.\ \par
\end{Lemma}
{\hskip 1.8em}Proof\ \par
for a sequence $S$ we introduce a related sequence $\{\epsilon _{a}\}_{a\in S}$ 
of independent random variables with the same law $P(\epsilon _{a}=1)=P(\epsilon _{a}=-1)=1/2$. 
We shall denote by ${\mathbb{E}}$ the associated expectation.\ \par
{\hskip 1.8em}Let $S$ be a $q$-Carleson sequence, with the associated 
$\{\epsilon _{a}\}_{a\in S}$ we have\ \par
{\hskip 3.6em}$\displaystyle \displaystyle \left\Vert{\displaystyle \sum_{a\in S}^{}{\mu _{a}\epsilon 
_{a}k_{q,a}}}\right\Vert _{q}^{q}\lesssim \displaystyle \left\Vert{\mu 
}\right\Vert _{q}^{q}$\ \par
because $\left\vert{\epsilon _{a}}\right\vert =1$. Taking expectation 
on both sides leads to\ \par
{\hskip 3.6em}$\displaystyle \displaystyle \left\Vert{{\mathbb{E}}\displaystyle \left[{\displaystyle 
\left\vert{\displaystyle \sum_{a\in S}^{}{\mu _{a}\epsilon _{a}k_{q,a}}}\right\vert 
^{q}}\right] }\right\Vert _{1}={\mathbb{E}}\displaystyle \left[{\displaystyle 
\left\Vert{\displaystyle \sum_{a\in S}^{}{\mu _{a}\epsilon _{a}k_{q,a}}}\right\Vert 
_{q}^{q}}\right] \lesssim \displaystyle \left\Vert{\mu }\right\Vert _{q}^{q}.$\ 
\par
Now using Khintchine's inequalities for the left expression\ \par
{\hskip 3.6em}$\displaystyle \displaystyle \left\Vert{{\mathbb{E}}\displaystyle \left[{\displaystyle 
\left\vert{\displaystyle \sum_{a\in S}^{}{\mu _{a}\epsilon _{a}k_{q,a}}}\right\vert 
^{q}}\right] }\right\Vert _{1}\simeq \displaystyle \left\Vert{\displaystyle 
\sum_{a\in S}^{}{\displaystyle \left\vert{\mu _{a}}\right\vert ^{2}\displaystyle 
\left\vert{k_{q,a}}\right\vert ^{2}}}\right\Vert _{q/2}^{q/2},$\ \par
we get\ \par
{\hskip 3.6em}$\displaystyle \displaystyle \left\Vert{\displaystyle \sum_{a\in S}^{}{\displaystyle 
\left\vert{\mu _{a}}\right\vert ^{2}\displaystyle \left\vert{k_{q,a}}\right\vert 
^{2}}}\right\Vert _{q/2}^{q/2}\lesssim {\mathbb{E}}\displaystyle \left[{\displaystyle 
\left\Vert{\displaystyle \sum_{a\in S}^{}{\mu _{a}\epsilon _{a}k_{q,a}}}\right\Vert 
_{q}^{q}}\right] \lesssim \displaystyle \left\Vert{\mu }\right\Vert _{q}^{q},$\ 
\par
and the lemma. \hfill$\square$\ \par
\ \par
{\hskip 1.8em}Now if $S$ is weakly $p$-Carleson is $S$ weakly $q$-Carleson 
for other $q$ ?\ \par
\ \par
Notice that any sequence $S$ is weakly $2$-Carleson :\ \par
{\hskip 3.6em}$\displaystyle \forall \nu \in \ell ^{1},\ \displaystyle \left\Vert{\displaystyle \sum_{a\in 
S}^{}{\nu _{a}\displaystyle \left\vert{k_{2,a}}\right\vert ^{2}}}\right\Vert 
_{1}\leq \displaystyle \sum_{a\in S}^{}{\displaystyle \left\vert{\nu _{a}}\right\vert 
\displaystyle \left\Vert{\displaystyle \left\vert{k_{2,a}}\right\vert 
^{2}}\right\Vert _{1}\leq \displaystyle \left\Vert{\nu }\right\Vert _{1}},$\ 
\par
because $\displaystyle \displaystyle \left\Vert{k_{2,a}}\right\Vert _{2}=\displaystyle \left\Vert{\displaystyle 
\left\vert{k_{2,a}}\right\vert ^{2}}\right\Vert _{1}=1$.\ \par
{\hskip 1.8em}Hence if $S$ is weakly $q$-Carleson with $q>2$ we can try 
to use interpolation of linear operators.\ \par
{\hskip 1.8em}Let us define our operator $T$ :\ \par
{\hskip 3.6em}$\displaystyle T\ :\ \ell ^{q}(\omega _{q}){\longrightarrow}L^{q}(\sigma );\ T\lambda 
:=\displaystyle \sum_{a\in S}^{}{\lambda _{a}\displaystyle \left\vert{k_{a}}\right\vert 
^{2}},$\ \par
with the weight $\omega _{q}(a):=\left\Vert{k_{a}}\right\Vert _{2q}^{-2q}$; 
this means that\ \par
{\hskip 2.1em}$\displaystyle \lambda \in \ell ^{q}(\omega _{q}){\Longrightarrow}\displaystyle \left\Vert{\lambda 
}\right\Vert _{\ell ^{q}(\omega _{q})}^{q}:=\displaystyle \sum_{a\in S}^{}{\displaystyle 
\left\vert{\lambda _{a}}\right\vert ^{q}\omega _{q}(a)}<\infty .$\ \par
{\hskip 1.8em}By a theorem of E. Stein and G. Weiss~\cite{SteinWeiss58} 
we know that if $T$ is bounded from $\displaystyle \ell ^{q}(\omega _{q})$ 
to $\displaystyle L^{q}(\sigma )$ and from $\displaystyle \ell ^{1}(\omega _{1})$ 
to $\displaystyle L^{1}(\sigma )$ then $T$ is bounded from $\displaystyle 
\ell ^{p}(\omega '_{p})$ to $\displaystyle L^{p}(\sigma )$ with $1\leq p\leq q$ 
provided that the weight satisfies the condition\ \par
{\hskip 3.6em}if $\displaystyle \displaystyle \frac{1}{p}=\displaystyle \frac{1-\theta }{1}+\displaystyle 
\frac{\theta }{q}$ then $\displaystyle \omega '_{p}=\omega _{1}^{p(1-\theta )}\omega _{q}^{p\theta /q}$.\ 
\par
{\hskip 1.8em}Here this means\ \par
{\hskip 3.6em}$\displaystyle \omega '_{p}(a)=\displaystyle \left\Vert{k_{a}}\right\Vert _{2}^{-2p(1-\theta 
)}\displaystyle \left\Vert{k_{a}}\right\Vert _{2q}^{-2p\theta }.$\ \par
Then $\displaystyle \displaystyle \left\Vert{T\lambda }\right\Vert _{p}^{p}\lesssim \displaystyle 
\left\Vert{\lambda }\right\Vert _{\ell ^{p}(\omega '_{p})}^{p}=\displaystyle 
\sum_{a\in S}^{}{\displaystyle \left\vert{\lambda _{a}}\right\vert ^{p}\omega 
'_{p}(a)}$. Hence if $\omega '_{p}(a)\lesssim \omega _{p}(a)$ we shall 
have\ \par
{\hskip 3.6em}$\displaystyle \displaystyle \left\Vert{T\lambda }\right\Vert _{p}^{p}\lesssim \displaystyle 
\left\Vert{\lambda }\right\Vert _{\ell ^{p}(\omega '_{p})}^{p}=\displaystyle 
\sum_{a\in S}^{}{\displaystyle \left\vert{\lambda _{a}}\right\vert ^{p}\omega 
'_{p}(a)}\lesssim \displaystyle \sum_{a\in S}^{}{\displaystyle \left\vert{\lambda 
_{a}}\right\vert ^{p}\omega _{p}(a)},$\ \par
and this will be OK.\ \par
\begin{Lemma} Let $q\geq 1$ and $\displaystyle \displaystyle \frac{1}{p}=\displaystyle \frac{1-\theta }{1}+\displaystyle 
\frac{\theta }{q}$ with $0<\theta <1$, then\ \par
{\hskip 3.6em}$\displaystyle \displaystyle \left\Vert{k_{a}}\right\Vert _{2p}\leq \displaystyle \left\Vert{k_{a}}\right\Vert 
_{2}^{(1-\theta )}\displaystyle \left\Vert{k_{a}}\right\Vert _{2q}^{\theta 
},$\ \par
\end{Lemma}
{\hskip 1.8em}Proof\ \par
let $\displaystyle \displaystyle \frac{1}{p}=\displaystyle \frac{1-\theta }{1}+\displaystyle 
\frac{\theta }{q}=\displaystyle \frac{1}{s}+\displaystyle \frac{1}{r}$ 
with $\displaystyle s=\displaystyle \frac{1}{1-\theta }$ and $\displaystyle 
r=\displaystyle \frac{q}{\theta }$.\ \par
{\hskip 1.8em}H{\"o}lder's inequalities give, for $f\in L^{s}(\sigma ),\ g\in L^{r}(\sigma )$\ 
\par
{\hskip 3.6em}$\displaystyle \displaystyle \left({\displaystyle \int_{X}^{}{\displaystyle \left\vert{fg}\right\vert 
^{p}\,d\sigma }}\right) ^{1/p}\leq \displaystyle \left({\displaystyle 
\int_{X}^{}{\displaystyle \left\vert{f}\right\vert ^{s}\,d\sigma }}\right) 
^{1/s}\displaystyle \left({\displaystyle \int_{X}^{}{\displaystyle \left\vert{g}\right\vert 
^{r}\,d\sigma }}\right) ^{1/r}.$\ \par
Set $f=\left\vert{k_{a}}\right\vert ^{2(1-\theta )},\ \ g:=\left\vert{k_{a}}\right\vert 
^{2\theta }$ we get\ \par
{\hskip 2.1em}$\displaystyle \displaystyle \left({\displaystyle \int_{X}^{}{\displaystyle \left\vert{k_{a}}\right\vert 
^{2p}\,d\sigma }}\right) ^{1/p}\leq \displaystyle \left({\displaystyle 
\int_{X}^{}{\displaystyle \left\vert{k_{a}}\right\vert ^{2(1-\theta )s}\,d\sigma 
}}\right) ^{1/s}\displaystyle \left({\displaystyle \int_{X}^{}{\displaystyle 
\left\vert{k_{a}}\right\vert ^{2\theta r}\,d\sigma }}\right) ^{1/r},$\ 
\par
hence replacing $s,r$\ \par
{\hskip 1.8em}$\displaystyle \displaystyle \left({\displaystyle \int_{X}^{}{\displaystyle \left\vert{k_{a}}\right\vert 
^{2p}\,d\sigma }}\right) ^{1/p}\leq \displaystyle \left({\displaystyle 
\int_{X}^{}{\displaystyle \left\vert{k_{a}}\right\vert ^{2}\,d\sigma }}\right) 
^{1-\theta }\displaystyle \left({\displaystyle \int_{X}^{}{\displaystyle 
\left\vert{k_{a}}\right\vert ^{2q}\,d\sigma }}\right) ^{\theta /q},$\ 
\par
hence\ \par
{\hskip 1.8em}$\displaystyle \displaystyle \left\Vert{k_{a}}\right\Vert _{2p}\leq \displaystyle \left\Vert{k_{a}}\right\Vert 
_{2}^{(1-\theta )}\displaystyle \left\Vert{k_{a}}\right\Vert _{2q}^{\theta 
},$\ \par
and the lemma. \hfill$\square$\ \par
{\hskip 1.8em}Back to our operator $T$, we have $\displaystyle \omega '_{p}(a)=\displaystyle \left\Vert{k_{a}}\right\Vert _{2}^{-2p(1-\theta 
)}\displaystyle \left\Vert{k_{a}}\right\Vert _{2q}^{-2p\theta }$ but the 
lemma above says $\displaystyle \displaystyle \left\Vert{k_{a}}\right\Vert _{2p}\lesssim \displaystyle 
\left\Vert{k_{a}}\right\Vert _{2}^{1-\theta }\displaystyle \left\Vert{k_{a}}\right\Vert 
_{2q}^{\theta }$ with $\displaystyle \displaystyle \frac{1}{p}=\displaystyle \frac{1-\theta }{1}+\displaystyle 
\frac{\theta }{q}$ which implies $\omega '_{p}(a)\lesssim \left\Vert{k_{a}}\right\Vert _{2p}^{-2p}=\omega 
_{p}(a)$ and the condition of the Stein-Weiss theorem are fullfilled, 
so we proved\ \par
\begin{Lemma} \label{algUniExt719}If $S$ is weakly $q$-Carleson, with 
$q>2$ then $S$ is weakly $p$-Carleson for any $2\leq p\leq q$.\ \par
\end{Lemma}
{\hskip 1.8em}We notice too that any sequence $S$ is $1$-Carleson\ \par
{\hskip 3.6em}$\displaystyle \forall \mu \in \ell ^{1},\ \displaystyle \left\Vert{\displaystyle \sum_{a\in 
S}^{}{\mu _{a}k_{1,a}}}\right\Vert _{1}\leq \displaystyle \sum_{a\in S}^{}{\displaystyle 
\left\vert{\mu _{a}}\right\vert \displaystyle \left\Vert{k_{1,a}}\right\Vert 
_{1}\leq \displaystyle \left\Vert{\mu }\right\Vert _{1}},$\ \par
{\hskip 1.8em}and the same proof as above gives\ \par
\begin{Lemma} If $S$ is $q$-Carleson, with $q>1$ then $S$ is $p$-Carleson 
for any $1\leq p\leq q$.\ \par
\end{Lemma}
\ \par
{\hskip 1.8em}In the ball or in the polydisc, we have much better:\ \par
\begin{Remark}  If $S$ is $q$-Carleson for a $q\in ]1,\infty [$ then $S$ 
is $p$-Carleson for any $p$. Moreover $S\ q$-Carleson is equivalent to 
$S$ weakly $2q$-Carleson.\ \par
\end{Remark}
\section{Main results{\hskip 1.8em}}
\setcounter{equation}{0}Now we are in position to state our main results.\ 
\par
\begin{Theorem} \label{algUniExt26}Let $p\geq 1$, $1\leq s<p$ and $q$ 
be such that $\displaystyle \displaystyle \frac{1}{s}=\displaystyle \frac{1}{p}+\displaystyle \frac{1}{q}$. 
Suppose that $S\subset {\mathcal{M}}_{s}\cap {\mathcal{M}}_{q'}$, that 
$S$ is dual bounded in $H^{p}(\sigma ),\ p\leq 2$., that $S$ is weakly 
$q$-Carleson and $\sigma $ verifies the structural hypotheses $SH(q)$ 
and $SH(p,s)$. Then $S$ is $\displaystyle H^{s}(\sigma )$ interpolating 
and has the $L.E.P.$ in $H^{s}(\sigma )$.\ \par
\end{Theorem}
{\hskip 1.8em}Using this time the fact that Kinchine's inequalities also 
provide a way to put absolut values inside sums, we get the other extremity 
for the range of $p$'s:\ \par
\begin{Theorem} \label{algUniExt27}Let $1\leq s<\infty $. Suppose that 
$S\subset {\mathcal{M}}_{s}\cap {\mathcal{M}}_{s'}$, that $S$ is dual 
bounded in $H^{\infty }(\sigma )$, $S$ is weakly $p$-Carleson for a $p>s$ 
and $(A,\ \sigma )$ verify the structural hypotheses $SH()$. Then $S$ 
is $\displaystyle H^{s}(\sigma )$ interpolating with the $L.E.P.$.\ \par
\end{Theorem}
{\hskip 1.8em}These theorems will be consequence of the next lemma.\ \par
{\hskip 1.5em}As above, if $S$ is a sequence of points in ${\mathcal{M}}$, 
we introduce the related sequence $\{\epsilon _{a}\}_{a\in S}$ of independent 
Bernouilli variables.\ \par
\begin{Lemma} \label{algUniExt411}Let $S\subset {\mathcal{M}}_{p}$ be 
a sequence of points such that a dual system $\{\rho _{p,a}\}_{a\in S}$ 
exists in $\displaystyle H^{p}(\sigma )$; let $1\leq s<p$ and $q$ be such 
that $\displaystyle \displaystyle \frac{1}{s}=\displaystyle \frac{1}{p}+\displaystyle \frac{1}{q}$;\ 
\par
if  $\displaystyle \forall \lambda \in \ell ^{p}(S),\ {\mathbb{E}}\displaystyle \left[{\displaystyle 
\left\Vert{\displaystyle \sum_{a\in S}^{}{\lambda _{a}\epsilon _{a}\rho 
_{p,a}}}\right\Vert _{p}^{p}}\right] \lesssim \displaystyle \left\Vert{\lambda 
}\right\Vert _{\ell ^{p}}^{p}$, $S$ is $q$-weakly Carleson and $\sigma $ 
verifies $SH(q)$, $SH(p,s)$ then $S$ is $\displaystyle H^{s}(\sigma )$ 
interpolating and moreover $S$ has the $L.E.P.$.\ \par
\end{Lemma}
{\hskip 1.8em}Proof\ \par
{\hskip 1.8em}If $p=1$ we have nothing to prove: the functions $\rho _{1,a}$ 
are uniformly bounded in $\displaystyle H^{1}(\sigma )$, just set\ \par
{\hskip 1.8em}$\displaystyle \forall \lambda \in \ell ^{1},\ T(\lambda ):=\displaystyle \sum_{a\in 
S}^{}{\lambda _{a}\rho _{1,a}},$\ \par
this function interpolates the sequence $\lambda $, is bounded in $\displaystyle 
H^{1}(\sigma )$, and clearly the operator $T$ is also linear and bounded.\ 
\par
{\hskip 1.8em}If $p>1$, we may suppose that $1<s<p$ because if $S\in IH^{s}(\sigma )$ 
then by theorem~\ref{algUniInt5}, for $S\subset {\mathcal{M}}_{1}$ we 
also have that $S\in IH^{1}(\sigma )$.\ \par
{\hskip 1.8em}First we truncate the sequence: $S_{N}$ is the first $N$ 
elements of $S$. We shall get estimates independent of $N$, i.e.\ \par
for $s\in [1,p[$ and $\nu \in \ell _{N}^{s}$ we shall built a function 
$h\in H^{s}(\sigma )$ such that:\ \par
{\hskip 1.5em}$\forall j=0,...,N-1,\ h(a_{j})=\nu _{j}\left\Vert{k_{a_{j}}}\right\Vert 
_{s'}$ and $\left\Vert{h}\right\Vert _{H^{s}}\leq C\left\Vert{\nu }\right\Vert _{\ell 
_{N}^{s}}$,\ \par
with the constant $C$ independent of $N$. We conclude then by use of lemma~\ref{algUniInt0}.\ 
\par
{\hskip 1.8em}We choose $q$ such that $\displaystyle \displaystyle \frac{1}{s}=\displaystyle \frac{1}{p}+\displaystyle \frac{1}{q}$; 
then $q\in ]p',\infty [$ with $p'$ the conjugate exponent of $p$ and we 
set $\nu _{j}=\lambda _{j}\mu _{j}$ with $\mu _{j}:=\left\vert{\nu _{j}}\right\vert ^{s/q}\in \ell ^{q},\ \lambda 
_{j}:=\frac{\nu _{j}}{\left\vert{\nu _{j}}\right\vert }\left\vert{\nu 
_{j}}\right\vert ^{s/p}\in \ell ^{p}$ then $\left\Vert{\nu }\right\Vert _{s}=\left\Vert{\lambda }\right\Vert _{p}\left\Vert{\mu 
}\right\Vert _{q}$.\ \par
Let $\displaystyle c_{a}:=\displaystyle \frac{\displaystyle \left\Vert{k_{a}}\right\Vert 
_{s'}}{\displaystyle \left\Vert{k_{a}}\right\Vert _{p'}k_{q,a}(a)}=\displaystyle 
\frac{\displaystyle \left\Vert{k_{a}}\right\Vert _{s'}\displaystyle \left\Vert{k_{a}}\right\Vert 
_{q}}{\displaystyle \left\Vert{k_{a}}\right\Vert _{p'}k_{a}(a)}$. By $SH(q)$ 
we have $k_{a}(a)\geq \alpha \left\Vert{k_{a}}\right\Vert _{q}\left\Vert{k_{a}}\right\Vert 
_{q'}$ hence\ \par
{\hskip 2.1em}$c_{a}\leq \frac{\left\Vert{k_{a}}\right\Vert _{s'}}{\alpha ^{-1}\left\Vert{k_{a}}\right\Vert 
_{p'}\left\Vert{k_{a}}\right\Vert _{q'}}$ and by $SH(p,s)$ we get $c_{a}\leq \alpha ^{-1}\beta $.\ 
\par
(i)  Now set $h(z)=\sum_{a\in S}^{}{\nu _{a}c_{a}\rho _{a}k_{q,a}}$ then:\ 
\par
{\hskip 3.6em}$\displaystyle h(a)=\nu _{a}\displaystyle \left\Vert{k_{a}}\right\Vert _{s'}$ 
because $\rho _{a}(b)=\delta _{ab}\left\Vert{k_{a}}\right\Vert _{p'}$.\ 
\par
{\hskip 2.1em}These are the good values, hence $\displaystyle h$ interpolates 
$\nu $ and moreover $h$ is clearly linear in $\nu $.\ \par
(ii)  Estimate on the $H^{s}(\sigma )$ norm of $h$.\ \par
{\hskip 1.8em}Set\ \par
{\hskip 2.1em}$f(\epsilon ,z):=\sum_{a\in S}^{}{\lambda _{a}c_{a}\epsilon _{a}\rho _{a}(z)},\ 
\ \ \ \ \ \ \ \ \ \ \ g(\epsilon ,z):=\sum_{a\in S}^{}{\mu _{a}\epsilon 
_{a}k_{q,a}(z)}.$\ \par
Then $\displaystyle h(z)={\mathbb{E}}(f(\epsilon ,z)g(\epsilon ,z))$ because 
$\displaystyle {\mathbb{E}}(\epsilon _{j}\epsilon _{k})=\delta _{jk}$.\ 
\par
So we get\ \par
{\hskip 1.8em}$\displaystyle \displaystyle \left\vert{h(z)}\right\vert ^{s}=\displaystyle \left\vert{{\mathbb{E}}(fg)}\right\vert 
^{s}\leq ({\mathbb{E}}(\displaystyle \left\vert{fg}\right\vert ))^{s}\leq 
{\mathbb{E}}(\displaystyle \left\vert{fg}\right\vert ^{s}),$\ \par
hence\ \par
{\hskip 1.8em}$\displaystyle \displaystyle \left\Vert{h}\right\Vert _{s}=\displaystyle \left({\displaystyle 
\int_{X}^{}{\displaystyle \left\vert{h(z)}\right\vert ^{s}\,d\sigma (z)}}\right) 
^{1/s}\leq \displaystyle \left({\displaystyle \int_{X}^{}{{\mathbb{E}}(\displaystyle 
\left\vert{fg}\right\vert ^{s})\,d\sigma (z)}}\right) ^{1/s}.$\ \par
But, using H{\"o}lder's inequality, we get\ \par
{\hskip 1.8em}
\begin{equation} 
\displaystyle \int_{X}^{}{{\mathbb{E}}(\displaystyle \left\vert{fg}\right\vert 
^{s})\,d\sigma (z)}={\mathbb{E}}\displaystyle \left[{\displaystyle \int_{X}^{}{\displaystyle 
\left\vert{fg}\right\vert ^{s}\,d\sigma (z)}}\right] \leq \displaystyle 
\left({{\mathbb{E}}\displaystyle \left[{\displaystyle \int_{X}^{}{\displaystyle 
\left\vert{f}\right\vert ^{p}\,d\sigma }}\right] }\right) ^{s/p}\displaystyle 
\left({{\mathbb{E}}\displaystyle \left[{\displaystyle \int_{X}^{}{\displaystyle 
\left\vert{g}\right\vert ^{q}\,d\sigma }}\right] }\right) ^{s/q}.\label{algUniExt410}
\end{equation} \ \par
{\hskip 1.8em}Let $\forall a\in S,\ \tilde \lambda _{a}:=c_{a}\lambda _{a}{\Longrightarrow}\left\Vert{\tilde 
\lambda }\right\Vert _{p}\leq \alpha \beta \left\Vert{\lambda }\right\Vert 
_{p}$ and the first factor is controlled by the lemma hypothesis\ \par
{\hskip 3.6em}
\begin{equation} 
{\mathbb{E}}\displaystyle \left[{\displaystyle \int_{X}^{}{\displaystyle 
\left\vert{f}\right\vert ^{p}\,d\sigma }}\right] ={\mathbb{E}}\displaystyle 
\left[{\displaystyle \left\Vert{\displaystyle \sum_{a\in S}^{}{\lambda 
_{a}c_{a}\epsilon _{a}\rho _{p,a}}}\right\Vert _{p}^{p}}\right] \lesssim 
\displaystyle \left\Vert{\tilde \lambda }\right\Vert _{p}^{p}\lesssim 
\displaystyle \left\Vert{\lambda }\right\Vert _{\ell ^{p}}^{p}.\label{algUniExt48}
\end{equation} \ \par
{\hskip 1.8em}Fubini theorem gives for the second factor\ \par
{\hskip 3.9em}$\displaystyle {\mathbb{E}}\displaystyle \left[{\displaystyle \int_{X}^{}{\displaystyle 
\left\vert{g}\right\vert ^{q}\,d\sigma }}\right] =\displaystyle \int_{X}^{}{{\mathbb{E}}\displaystyle 
\left[{\displaystyle \left\vert{g}\right\vert ^{q}}\right] \,d\sigma }.$\ 
\par
{\hskip 1.8em}We apply Khintchine's inequalities to $\displaystyle {\mathbb{E}}\displaystyle \left[{\displaystyle \left\vert{g}\right\vert 
^{q}}\right] $\ \par
{\hskip 3.6em}$\displaystyle {\mathbb{E}}\displaystyle \left[{\displaystyle \left\vert{g}\right\vert 
^{q}}\right] \simeq \displaystyle \left({\displaystyle \sum_{a\in S}^{}{\displaystyle 
\left\vert{\mu _{a}}\right\vert ^{2}\displaystyle \left\vert{k_{q,a}}\right\vert 
^{2}}}\right) ^{q/2},$\ \par
hence $S$ being weak $q$-Carleson implies\ \par
{\hskip 3.6em}
\begin{equation} 
\displaystyle \int_{X}^{}{{\mathbb{E}}\displaystyle \left[{\displaystyle 
\left\vert{g}\right\vert ^{q}}\right] \,d\sigma }\lesssim \displaystyle 
\int_{X}^{}{\displaystyle \left({\displaystyle \sum_{a\in S}^{}{\displaystyle 
\left\vert{\mu _{a}}\right\vert ^{2}\displaystyle \left\vert{k_{q,a}}\right\vert 
^{2}}}\right) ^{q/2}\,d\sigma }\lesssim \displaystyle \left\Vert{\mu }\right\Vert 
_{\ell ^{q}}^{q}.\label{algUniExt49}
\end{equation} \ \par
{\hskip 1.8em}So putting~(\ref{algUniExt48}) and~(\ref{algUniExt49}) in~(\ref{algUniExt410}) 
we get the lemma. \hfill$\square$\ \par
\subsection{Proof of theorem~\ref{algUniExt26}.{\hskip 1.8em}}
Let us recall the theorem we want to prove.\ \par
\begin{Theorem} Let $p\geq 1$, $1\leq s<p$ and $q$ be such that $\displaystyle 
\displaystyle \frac{1}{s}=\displaystyle \frac{1}{p}+\displaystyle \frac{1}{q}$. 
Suppose that $S\subset {\mathcal{M}}_{s}\cap {\mathcal{M}}_{q'}$, that 
$\{\rho _{p,a}\}_{a\in S}$ is a norm bounded sequence in $H^{p}(\sigma ),\ p\leq 2$, 
that $S$ is weakly $q$-Carleson and $\sigma $ verifies the structural 
hypotheses $SH(q),\ SH(p,s)$. Then $S$ is $H^{s}(\sigma )$-interpolating 
with the $L.E.P.$.\ \par
\end{Theorem}
{\hskip 1.8em}It remains to prove that the hypotheses of the theorem implies 
those of the lemma~\ref{algUniExt411}.\ \par
{\hskip 1.8em}We have to prove that\ \par
{\hskip 3.6em}$\displaystyle {\mathbb{E}}\displaystyle \left[{\displaystyle \left\Vert{\displaystyle 
\sum_{a\in S}^{}{\lambda _{a}\epsilon _{a}\rho _{p,a}}}\right\Vert _{p}^{p}}\right] 
\lesssim \displaystyle \left\Vert{\lambda }\right\Vert _{\ell ^{p}}^{p},$\ 
\par
knowing that the dual sequence $\{\rho _{p,a}\}_{a\in S}$ is bounded in 
$\displaystyle H^{p}(\sigma )$, i.e.\ \par
{\hskip 3.6em}$\displaystyle \sup _{a\in S}\ \displaystyle \left\Vert{\rho _{p,a}}\right\Vert _{p}\leq 
C.$\ \par
By Fubini's theorem\ \par
{\hskip 3.9em}$\displaystyle {\mathbb{E}}\displaystyle \left[{\displaystyle \left\Vert{\displaystyle 
\sum_{a\in S}^{}{\lambda _{a}\epsilon _{a}\rho _{p,a}}}\right\Vert _{p}^{p}}\right] 
=\displaystyle \int_{X}^{}{{\mathbb{E}}\displaystyle \left[{\displaystyle 
\left\vert{\displaystyle \sum_{a\in S}^{}{\lambda _{a}\epsilon _{a}\rho 
_{p,a}}}\right\vert ^{p}}\right] \,d\sigma },$\ \par
and by Khintchine's inequalities we have\ \par
{\hskip 3.6em}$\displaystyle {\mathbb{E}}\displaystyle \left[{\displaystyle \left\vert{\displaystyle 
\sum_{a\in S}^{}{\lambda _{a}\epsilon _{a}\rho _{p,a}}}\right\vert ^{p}}\right] 
\simeq \displaystyle \left({\displaystyle \sum_{a\in S}^{}{\displaystyle 
\left\vert{\lambda _{a}}\right\vert ^{2}\displaystyle \left\vert{\rho 
_{p,a}}\right\vert ^{2}}}\right) ^{p/2}.$\ \par
Now $p\leq 2$, so $\displaystyle \displaystyle \left({\displaystyle \sum_{a\in S}^{}{\displaystyle \left\vert{\lambda 
_{a}}\right\vert ^{2}\displaystyle \left\vert{\rho _{p,a}}\right\vert 
^{2}}}\right) ^{1/2}\leq \displaystyle \left({\displaystyle \sum_{a\in 
S}^{}{\displaystyle \left\vert{\lambda _{a}}\right\vert ^{p}\displaystyle 
\left\vert{\rho _{p,a}}\right\vert ^{p}}}\right) ^{1/p}$ hence\ \par
{\hskip 1.2em}$\displaystyle \displaystyle \int_{X}^{}{{\mathbb{E}}\displaystyle \left[{\displaystyle 
\left\vert{\displaystyle \sum_{a\in S}^{}{\lambda _{a}\epsilon _{a}\rho 
_{p,a}}}\right\vert ^{p}}\right] \,d\sigma }\leq \displaystyle \int_{X}^{}{\displaystyle 
\left({\displaystyle \sum_{a\in S}^{}{\displaystyle \left\vert{\lambda 
_{a}}\right\vert ^{p}\displaystyle \left\vert{\rho _{p,a}}\right\vert 
^{p}}}\right) }\,d\sigma =\displaystyle \sum_{a\in S}^{}{\displaystyle 
\left\vert{\lambda _{a}}\right\vert ^{p}\displaystyle \left\Vert{\rho 
_{p,a}}\right\Vert _{p}^{p}}.$\ \par
{\hskip 1.8em}So, finally\ \par
{\hskip 3.6em}$\displaystyle {\mathbb{E}}\displaystyle \left[{\displaystyle \left\Vert{\displaystyle 
\sum_{a\in S}^{}{\lambda _{a}\epsilon _{a}\rho _{p,a}}}\right\Vert _{p}^{p}}\right] 
\lesssim \sup _{a\in S}\ \displaystyle \left\Vert{\rho _{p,a}}\right\Vert 
_{p}^{p}\displaystyle \left\Vert{\lambda }\right\Vert _{p}^{p},$\ \par
and the theorem~\ref{algUniExt26}. \hfill$\square$\ \par
{\hskip 1.8em}Suggested by F. Lust-Piquard, one can use that $H^{p}(\sigma )\subset L^{p}(\sigma )$ 
hence, because $p\leq 2$, $\displaystyle H^{p}(\sigma )$ is of type $p$ 
which means precisely(~\cite{LiQuef04}, Th III.9) that $\displaystyle 
{\mathbb{E}}\displaystyle \left[{\displaystyle \left\vert{\displaystyle 
\sum_{a\in S}^{}{\lambda _{a}\epsilon _{a}\rho _{p,a}}}\right\vert ^{p}}\right] 
\lesssim \displaystyle \sum_{a\in S}^{}{\displaystyle \left\vert{\lambda 
_{a}\rho _{a}}\right\vert ^{p}},$ hence integrating and using Fubini, 
we get\ \par

\begin{displaymath} 
{\mathbb{E}}\displaystyle \left[{\displaystyle \left\Vert{\displaystyle 
\sum_{a\in S}^{}{\lambda _{a}\epsilon _{a}\rho _{p,a}}}\right\Vert _{p}^{p}}\right] 
\lesssim \displaystyle \int_{}^{}{\displaystyle \sum_{a\in S}^{}{\displaystyle 
\left\vert{\lambda _{a}\rho _{a}}\right\vert ^{p}}\,d\sigma }\lesssim 
(\sup _{a\in S}\ \displaystyle \left\Vert{\rho _{p,a}}\right\Vert _{p})\displaystyle 
\left\Vert{\lambda }\right\Vert _{\ell ^{p}}^{p}\lesssim \displaystyle 
\left\Vert{\lambda }\right\Vert _{\ell ^{p}}^{p}.\end{displaymath} \ \par
And again the theorem.\ \par
\subsection{Proof of theorem~\ref{algUniExt27}.{\hskip 1.8em}}
Let us recall the theorem we want to prove.\ \par
\begin{Theorem} Let $1\leq s<\infty $. Suppose that $S\subset {\mathcal{M}}_{s}\cap {\mathcal{M}}_{s'}$, 
that $\{\rho _{a}\}_{a\in S}$  is a norm bounded sequence in $H^{\infty }(\sigma )$, 
weakly $p$-Carleson for a $p>s$ and $\sigma $ verifies the structural 
hypotheses $SH(p,s),\ SH(q)$ for $q$ such that $\displaystyle \displaystyle \frac{1}{s}=\displaystyle \frac{1}{p}+\displaystyle \frac{1}{q}$. 
Then $S$ is $\displaystyle H^{s}(\sigma )$ interpolating with the $L.E.P.$.\ 
\par
\end{Theorem}
{\hskip 1.5em}Proof\ \par
the idea is still to use lemma~\ref{algUniExt411}, but in two steps. Let 
$s<\infty $ be given and take $p$ such that $s<p<\infty $ and $S$ is weakly 
$p$-Carleson.\ \par
{\hskip 1.8em}Set $\displaystyle \forall a\in S,\ \rho _{p,a}:=\rho _{a}k_{p,a}$. 
We have $\displaystyle \displaystyle \left\Vert{\rho _{p,a}}\right\Vert _{p}\leq \displaystyle 
\left\Vert{\rho _{a}}\right\Vert _{\infty }\displaystyle \left\Vert{k_{p,a}}\right\Vert 
_{p}=\displaystyle \left\Vert{\rho _{a}}\right\Vert _{\infty }\leq C$ 
by hypothesis. We want to prove that\ \par
{\hskip 3.6em}$\displaystyle {\mathbb{E}}\displaystyle \left[{\displaystyle \left\Vert{\displaystyle 
\sum_{a\in S}^{}{\lambda _{a}\epsilon _{a}\rho _{p,a}}}\right\Vert _{p}^{p}}\right] 
=\displaystyle \int_{X}^{}{{\mathbb{E}}\displaystyle \left[{\displaystyle 
\left\vert{\displaystyle \sum_{a\in S}^{}{\lambda _{a}\epsilon _{a}\rho 
_{p,a}}}\right\vert ^{p}}\right] \,d\sigma }\lesssim \displaystyle \left\Vert{\lambda 
}\right\Vert _{\ell ^{p}}^{p},$\ \par
in order to apply lemma~\ref{algUniExt411}.\ \par
{\hskip 1.8em}By Khintchine's inequalities we have\ \par
{\hskip 1.8em}$\displaystyle {\mathbb{E}}\displaystyle \left[{\displaystyle \left\vert{\displaystyle 
\sum_{a\in S}^{}{\lambda _{a}\epsilon _{a}\rho _{p,a}}}\right\vert ^{p}}\right] 
\simeq \displaystyle \left({\displaystyle \sum_{a\in S}^{}{\displaystyle 
\left\vert{\lambda _{a}}\right\vert ^{2}\displaystyle \left\vert{\rho 
_{p,a}}\right\vert ^{2}}}\right) ^{p/2},$\ \par
but this time we use that $\displaystyle \displaystyle \left\vert{\rho _{p,a}}\right\vert \leq \displaystyle \left\Vert{\rho 
_{\infty ,a}}\right\Vert \displaystyle \left\vert{k_{a,p}}\right\vert 
\leq C\displaystyle \left\vert{k_{a,p}}\right\vert $ hence\ \par
{\hskip 3.6em}$\displaystyle {\mathbb{E}}\displaystyle \left[{\displaystyle \left\vert{\displaystyle 
\sum_{a\in S}^{}{\lambda _{a}\epsilon _{a}\rho _{p,a}}}\right\vert ^{p}}\right] 
\lesssim C^{p}\displaystyle \left({\displaystyle \sum_{a\in S}^{}{\displaystyle 
\left\vert{\lambda _{a}}\right\vert ^{2}\displaystyle \left\vert{k_{a,p}}\right\vert 
^{2}}}\right) ^{p/2}.$\ \par
Using that $S$ is weakly $p$-Carleson, we get\ \par
{\hskip 3.6em}$\displaystyle \displaystyle \left\Vert{\displaystyle \sum_{a\in S}^{}{\displaystyle 
\left\vert{\lambda _{a}}\right\vert ^{2}\displaystyle \left\vert{k_{a,p}}\right\vert 
^{2}}}\right\Vert _{p/2}^{p/2}\leq D\displaystyle \left\Vert{\lambda }\right\Vert 
_{p}^{p},$\ \par
hence\ \par
{\hskip 3.6em}$\displaystyle {\mathbb{E}}\displaystyle \left[{\displaystyle \left\Vert{\displaystyle 
\sum_{a\in S}^{}{\lambda _{a}\epsilon _{a}\rho _{p,a}}}\right\Vert _{p}^{p}}\right] 
\lesssim C^{p}\displaystyle \int_{X}^{}{\displaystyle \left({\displaystyle 
\sum_{a\in S}^{}{\displaystyle \left\vert{\lambda _{a}}\right\vert ^{2}\displaystyle 
\left\vert{k_{a,p}}\right\vert ^{2}}}\right) \,d\sigma }\lesssim DC^{p}\displaystyle 
\left\Vert{\lambda }\right\Vert _{\ell ^{p}}^{p},$\ \par
and we can apply the lemma~\ref{algUniExt411} which gives the theorem 
because $p>s$ implies that $S$ is still weakly $s$-Carleson by lemma~\ref{algUniExt719}. 
\hfill$\square$\ \par
\section{Application to the ball and to the polydisc.{\hskip 1.8em}}
\setcounter{equation}{0}In~\cite{amIntInt06} it is proved that the structural 
hypotheses hold in the polydisc. Moreover the Carleson measures, hence 
the Carleson sequences, are characterized geometrically and they are the 
same for all $p\in ]1,\infty [$ (see~\cite{AChang79},~\cite{ChangFeff85}). 
So it is enough to say "Carleson sequence" in the theorem:\ \par
\begin{Theorem} Let $S\subset {\mathbb{D}}^{n}$ be a Carleson sequence 
and dual bounded  in $\displaystyle H^{p}({\mathbb{D}}^{n})$ with either 
$p=\infty $ or $p\leq 2$, then $S$ is $\displaystyle H^{s}({\mathbb{D}}^{n})$ 
interpolating for any $s<p$ with the LEP.\ \par
\end{Theorem}
{\hskip 1.8em}Still in~\cite{amIntInt06} it is proved that the structural 
hypotheses hold in the ball. Again the Carleson measures, hence the Carleson 
sequences, are characterized geometrically and they are the same for all 
$p\in ]1,\infty [$ (see~\cite{HormPSH67})but moreover a theorem of P. 
Thomas~\cite{Thomas87} gives that $S$ dual bounded in $\displaystyle H^{p}({\mathbb{B}})$ 
implies $S$ Carleson, hence\ \par
\begin{Theorem} Let $S\subset {\mathbb{B}}$ be dual bounded  in $\displaystyle 
H^{p}({\mathbb{B}})$ with either $p=\infty $ or $p\leq 2$, then $S$ is 
$\displaystyle H^{s}({\mathbb{B}})$ interpolating for any $s<p$ with the 
LEP.\ \par
\end{Theorem}
{\hskip 1.8em}We have for free the same result for the Bergman classes 
of the ball by the "subordination lemma"~\cite{amBerg78}:\ \par
{\hskip 1.8em}to a function $f(z)$ defined on $z=(z_{1},...,z_{n})\in {\mathbb{B}}_{n}\subset {\mathbb{C}}^{n}$ 
associate the function.\ \par
{\hskip 3.6em}$\displaystyle \tilde f(z,w):=f(z)$ defined on $(z,w)=(z_{1},...,z_{n},w)\in {\mathbb{B}}_{n+1}\subset {\mathbb{C}}^{n+1}$.\ 
\par
Then we have that $f\in A^{p}({\mathbb{B}}_{n}){\Longleftrightarrow}\tilde f\in H^{p}({\mathbb{B}}_{n+1})$ 
with the same norm. Moreover if $F\in H^{p}({\mathbb{B}}_{n+1})$then $f(z):=F(z,0)\in A^{p}({\mathbb{B}}_{n})$ 
with $\displaystyle \displaystyle \left\Vert{f}\right\Vert _{A^{p}({\mathbb{B}}_{n})}\leq 
\displaystyle \left\Vert{F}\right\Vert _{H^{p}({\mathbb{B}}_{n+1})}$.\ 
\par
\ \par
Suppose that $S\subset {\mathbb{B}}_{n}$ is dual bounded in $\displaystyle 
A^{p}({\mathbb{B}}_{n})$ this means that\ \par
{\hskip 3.6em}$\displaystyle \exists \{\rho _{a}\}_{a\in S}\ s.t.\ \forall a\in S,\displaystyle \left\Vert{\rho 
_{a}}\right\Vert _{A^{p}({\mathbb{B}}_{n})}\leq C$ and $\displaystyle 
\rho _{a}(b)=\delta _{ab}(1-\displaystyle \left\vert{a}\right\vert ^{2})^{-(n+1)/p}$,\ 
\par
because the normalized reproducing kernel for $A^{p}({\mathbb{B}}_{n})$ 
is $\displaystyle b_{a}(z):=\displaystyle \frac{(1-\displaystyle \left\vert{a}\right\vert 
^{2})^{(n+1)/p'}}{(1-\overline{a}\cdot z)^{n+1}}$.\ \par
{\hskip 1.8em}Embed $S$ in ${\mathbb{B}}_{n+1}$ by $\tilde S:=\{(a,0),\ a\in S\}$ 
as in~\cite{amBerg78}, then the sequence $\{\tilde \rho _{a}\}_{a\in S}$ 
is precisely a bounded dual sequence for $\tilde S\subset {\mathbb{B}}_{n+1}$ 
in $\displaystyle H^{p}({\mathbb{B}}_{n+1})$ hence we can apply the previous 
theorem:\ \par
{\hskip 1.8em}if $p=\infty $ or $p\leq 2$ then $\tilde S$ is $H^{s}({\mathbb{B}}_{n+1})$ 
interpolating with the $L.E.P.$. If $T$ is the operator making the extension,\ 
\par
{\hskip 2.1em}$\forall \lambda \in \ell ^{s}{\longrightarrow}T\lambda \in H^{s}({\mathbb{B}}_{n+1}),\ 
(T\lambda )(a,0)=\lambda _{a}\left\Vert{k_{(a,0)}}\right\Vert _{H^{s'}({\mathbb{B}}_{n+1})},\ 
\left\Vert{T\lambda }\right\Vert _{H^{s}({\mathbb{B}}_{n+1})}\leq C_{I}\left\Vert{\lambda 
}\right\Vert _{s}$\ \par
then the operator $(U\lambda )(z):=(T\lambda )(z,0)$ is a bounded linear 
operator from $\ell ^{s}$ to $A^{s}({\mathbb{B}}_{n})$ making the extension 
because $\displaystyle \displaystyle \left\Vert{k_{(a,0)}}\right\Vert _{H^{s'}({\mathbb{B}}_{n+1})}=\displaystyle 
\left\Vert{b_{a}}\right\Vert _{A^{s'}({\mathbb{B}}_{b})}$ where $k$ is 
the kernel for $H^{s}({\mathbb{B}}_{n+1})$ and $b$ is the kernel for $A^{s}({\mathbb{B}}_{n})$. 
Hence we proved\ \par
\begin{Corollary} Let $S\subset {\mathbb{B}}$ be dual bounded  in $\displaystyle 
A^{p}({\mathbb{B}})$ with either $p=\infty $ or $p\leq 2$, then $S$ is 
$\displaystyle A^{s}({\mathbb{B}})$ interpolating for any $s<p$ with the 
LEP.\ \par
\end{Corollary}
{\hskip 1.8em}We also get the same result for the Bergman spaces with 
weight of the form $(1-\left\vert{z}\right\vert ^{2})^{k},\ k\in {\mathbb{N}}$ 
just by the same method, but considering $\displaystyle H^{p}({\mathbb{B}}_{n+k+1})$ 
instead of $\displaystyle H^{p}({\mathbb{B}}_{n+1})$.\ \par

\bibliographystyle{/usr/share/texmf/bibtex/bst/base/plain}

\end{document}